\documentclass[11pt]{article}
\usepackage[a4paper,top=2.7cm,bottom=2.9cm,left=2.9cm,right=2.9cm]{geometry}
\usepackage[T1]{fontenc}
\usepackage[utf8]{inputenc}
\usepackage[sc]{mathpazo}
\usepackage{amsmath,amssymb,amsthm,mathtools}
\usepackage{microtype}
\usepackage{enumitem}
\usepackage{titlesec}
\usepackage[hidelinks]{hyperref}

\titleformat{\section}{\large\bfseries}{\thesection}{0.7em}{}
\titleformat{\subsection}{\normalsize\bfseries}{\thesubsection}{0.6em}{}
\titlespacing*{\section}{0pt}{1.6\baselineskip}{0.7\baselineskip}
\titlespacing*{\subsection}{0pt}{1.2\baselineskip}{0.5\baselineskip}
\setlist{itemsep=2pt,topsep=4pt}
\allowdisplaybreaks
\newtheoremstyle{mainplain}{9pt}{9pt}{\itshape}{}{\bfseries}{.}{0.5em}{}
\newtheoremstyle{maindef}{9pt}{9pt}{}{}{\bfseries}{.}{0.5em}{}
\theoremstyle{mainplain}
\newtheorem{theorem}{Theorem}[section]
\newtheorem{proposition}[theorem]{Proposition}
\newtheorem{lemma}[theorem]{Lemma}
\newtheorem{corollary}[theorem]{Corollary}
\theoremstyle{maindef}
\newtheorem{definition}[theorem]{Definition}
\newtheorem{example}[theorem]{Example}
\newtheorem{remark}[theorem]{Remark}
\numberwithin{equation}{section}

\newcommand{\PP}{\mathbb{P}}
\newcommand{\LL}{\mathfrak{L}}
\newcommand{\Le}{\mathfrak{L}_\epsilon}
\newcommand{\M}{\mathcal{M}}
\newcommand{\Hq}{H_q}
\newcommand{\Kq}{\mathcal{K}_q}
\newcommand{\qn}[1]{[#1]_q}
\newcommand{\ip}[2]{\langle #1,#2\rangle}
\newcommand{\hqi}{h_{q^{-1}}}

\title{\bfseries A $q$-analogue of an inverse problem\\ leading to second-order linear functionals}
\author{Wathek Chammam\\[5pt]
\small Department of Mathematics, College of Science at Zulfi, Majmaah University,\\
\small Majmaah 11952, Saudi Arabia\\[2pt]
\small \href{mailto:w.chammam@mu.edu.sa}{w.chammam@mu.edu.sa}}
\date{}

\begin{document}
\maketitle

\begin{abstract}
\noindent
R.~Sfaxi studied the linear functional $\Le$ determined by the equation
$\Le-\epsilon\Le'=\LL$. We replace the derivative by the Hahn operator $\Hq$
and study the solution $\Le$ of $\Le-\epsilon\Hq\Le=\LL$. For $0<q<1$ we solve
the equation explicitly when $\LL$ is a Dirac mass, and we deduce that $\Le$
is positive-definite whenever $\epsilon>0$ and $\LL$ is positive-definite on
$[0,+\infty)$. We show that $\Le$ is $\Hq$-semiclassical if and only if
$E\LL=F\Le$ for two non-zero polynomials $E$ and $F$, and we work out the
$q$-analogues of the Hermite and Laguerre examples. Finally we treat the
$\Hq$-classical functional $\Kq$ whose moments are the $q$-factorials: here
$\Le=(1+\epsilon)^{-1}(\Kq+\epsilon\,h_{-\epsilon}\Kq)$, and for $\epsilon<0$
the functional $\Le$ is positive-definite if and only if $\epsilon=-q^j$ for
some non-zero integer $j$. It is $\Hq$-semiclassical for these values and for
no other non-zero real $\epsilon$. For $q=1$, by contrast,
positive-definiteness holds for every real $\epsilon$.
\end{abstract}

\medskip
\noindent\textbf{2020 Mathematics Subject Classification.} Primary 33D45; Secondary 42C05.

\smallskip
\noindent\textbf{Keywords.} Orthogonal polynomials; Hahn operator; semiclassical linear functionals;
positive-definite linear functionals; little $q$-Laguerre polynomials.

%======================================================================
\section{Introduction and preliminary results}\label{sec1}
%======================================================================

A quasi-definite linear functional on polynomials is called semiclassical
when it satisfies a first-order equation of Pearson type with polynomial
coefficients. In \cite{Sf}, Sfaxi started from a functional $\LL$ and
considered the functional $\Le$ defined by $\Le-\epsilon\Le'=\LL$. He
showed that $\Le$ is positive-definite as soon as $\LL$ is, gave a criterion
for $\Le$ to be semiclassical, and observed that, when $\LL$ is
semiclassical, $\Le$ satisfies a second-order equation which need not
reduce to a first-order one. The Laguerre case is treated there in detail.

In this note we ask what remains of this picture when the derivative is
replaced by the Hahn operator $\Hq$, in the setting of the
$\Hq$-semiclassical functionals \cite{Kh,KM}. Part of
it goes through with little change: the moments of $\Le$, the second-order
equation and the criterion $E\LL=F\Le$ are obtained much as in \cite{Sf},
and the formulas are $q$-versions of the original ones. Positivity is another
matter. The operator $\Hq$ does not commute with translations, and the sign
of $\epsilon$ now has to agree with the half-line that carries $\LL$. With
the other sign positivity may fail: in the example of Section~\ref{sec4} it
holds only when $\epsilon=-q^j$ for some non-zero integer $j$.

The paper is organized as follows. The rest of this section fixes the
notation and collects the rules of calculation for $\Hq$.
Section~\ref{sec2} contains the existence and uniqueness of $\Le$, the
explicit solution for a Dirac mass (Lemma~\ref{lem:delta}) and the
positivity theorem. Section~\ref{sec3} deals with the $\Hq$-semiclassical
character and with the $q$-Hermite and $q$-Laguerre examples. In
Section~\ref{sec4} the equation is solved completely for $\LL=\Kq$,
including the degenerate value $\epsilon=-1$.

Throughout, $\PP$ is the vector space of polynomials in one variable with
complex coefficients and $\PP'$ is its algebraic dual. We denote by
$\ip{\LL}{p}$ the action of $\LL\in\PP'$ on $p\in\PP$, and by
$(\LL)_n:=\ip{\LL}{x^n}$, $n\ge0$, the moments of $\LL$. A functional is
\emph{normalized} if $(\LL)_0=1$.

Let $q\in\mathbb C\setminus\{0\}$ be such that $q$ is not a root of unity.
We use the standard notation
\[
\qn n=\frac{1-q^n}{1-q},\qquad \qn0!=1,\qquad
\qn n!=\prod_{k=1}^n\qn k\quad(n\ge1),
\]
\[
(a;q)_0=1,\qquad (a;q)_n=\prod_{i=0}^{n-1}(1-aq^i),\qquad
(a;q)_\infty=\prod_{i\ge0}(1-aq^i)\quad(|q|<1).
\]
Since $q$ is not a root of unity, $\qn n\neq0$ for every $n\ge1$.

The Hahn operator $\Hq$ and the dilation $h_a$ ($a\in\mathbb C^*$) act on
$\PP$ by
\[
(\Hq f)(x)=\frac{f(qx)-f(x)}{(q-1)x},\qquad (h_af)(x)=f(ax).
\]
Thus $\Hq x^n=\qn n\,x^{n-1}$, so $\Hq$ lowers the degree by one, and
$\Hq f\to f'$ when $q\to1$. As usual \cite{KM,Sf}, these operators and the
multiplication by a polynomial $g$ are transposed to $\PP'$ by
\begin{equation}\label{eq:dual}
\ip{\Hq\LL}{f}:=-\ip{\LL}{\Hq f},\qquad
\ip{h_a\LL}{f}:=\ip{\LL}{h_af},\qquad
\ip{g\LL}{f}:=\ip{\LL}{gf}.
\end{equation}
In particular
\begin{equation}\label{eq:momH}
(\Hq\LL)_0=0,\qquad (\Hq\LL)_n=-\qn n\,(\LL)_{n-1}\quad(n\ge1),
\end{equation}
so that $\Hq\LL=0$ implies $\LL=0$. We also use the Dirac mass
$\ip{\delta_c}{f}=f(c)$.

A functional $\LL$ is \emph{quasi-definite} if there exists a monic
polynomial sequence $\{B_n\}_{n\ge0}$, $\deg B_n=n$, such that
$\ip{\LL}{B_nB_m}=r_n\delta_{n,m}$ with $r_n\neq0$ for all $n$; it is
\emph{positive-definite} if moreover $r_n>0$ for all $n$, which is
equivalent to $\ip{\LL}{p}>0$ for every non-zero real polynomial $p$ with
$p\ge0$ on $\mathbb R$ \cite{Ch}. We say that $\LL$ is
\emph{positive-definite on $[0,+\infty)$} if $\ip{\LL}{p}>0$ for every
non-zero real polynomial $p$ with $p\ge0$ on $[0,+\infty)$;
positive-definiteness on $(-\infty,0]$ is defined in the same way. Two
elementary facts will be used several times:
\begin{itemize}
\item if $\LL$ is quasi-definite and $A\LL=0$ for a polynomial $A$, then
$A=0$;
\item if $(x-c)\M=0$, then $\M=(\M)_0\,\delta_c$.
\end{itemize}

\begin{definition}[\cite{Kh,KM}]
A quasi-definite functional $\LL$ is \emph{$\Hq$-semiclassical} if there
are polynomials $\phi$, monic, and $\psi$, with $\deg\psi\ge1$, such that
\begin{equation}\label{eq:Pearson}
\Hq(\phi\LL)+\psi\LL=0 .
\end{equation}
It is \emph{$\Hq$-classical} if moreover $\deg\phi\le2$ and $\deg\psi=1$.
\end{definition}

By \eqref{eq:dual}, equation \eqref{eq:Pearson} means
$\ip{\LL}{\psi f-\phi\,\Hq f}=0$ for every $f\in\PP$.

\begin{lemma}\label{lem:rules}
Let $f,g\in\PP$, $\LL\in\PP'$ and $a\in\mathbb C^*$.
\begin{enumerate}[label=(\roman*)]
\item $\Hq(fg)=(h_qf)\,\Hq g+g\,\Hq f$.
\item $\Hq(f\LL)=(\hqi f)\,\Hq\LL+\Hq(\hqi f)\,\LL$; equivalently,
\begin{equation}\label{eq:fH}
f\,\Hq\LL=\Hq\big((h_qf)\LL\big)-(\Hq f)\,\LL .
\end{equation}
\item $h_a(\Hq\LL)=a\,\Hq(h_a\LL)$; in particular
$h_{-1}(\Hq\LL)=-\Hq(h_{-1}\LL)$.
\item Let $u\in\PP'$ satisfy $\Hq(\phi u)+\psi u=0$ and let $R\in\PP$. Then
$v=Ru$ satisfies
\[
\Hq\big((h_qR)\,\phi\,v\big)+\Big((\hqi R)\,\psi-\phi\,\Hq R
-\phi\,\Hq(\hqi R)\Big)v=0 .
\]
\end{enumerate}
\end{lemma}

\begin{proof}
(i) We have
\begin{align*}
(\Hq(fg))(x)&=\frac{f(qx)g(qx)-f(qx)g(x)}{(q-1)x}
+\frac{f(qx)g(x)-f(x)g(x)}{(q-1)x}\\
&=f(qx)\,(\Hq g)(x)+g(x)\,(\Hq f)(x).
\end{align*}

(ii) Let $g\in\PP$. Applying (i) to the polynomials $\hqi f$ and $g$,
and using $h_q\hqi f=f$, we get
$\Hq\big((\hqi f)g\big)=f\,\Hq g+g\,\Hq(\hqi f)$. Therefore
\[
\ip{(\hqi f)\Hq\LL}{g}=-\ip{\LL}{\Hq\big((\hqi f)g\big)}
=-\ip{\LL}{f\,\Hq g}-\ip{\LL}{g\,\Hq(\hqi f)}.
\]
The first term equals $\ip{\Hq(f\LL)}{g}$ and the second one equals
$-\ip{\Hq(\hqi f)\LL}{g}$, which proves the first identity. Replacing $f$
by $h_qf$ and using $\hqi h_qf=f$ gives \eqref{eq:fH}.

(iii) For $f\in\PP$, $\Hq\big(f(a\,\cdot)\big)(x)
=\frac{f(aqx)-f(ax)}{(q-1)x}=a\,(\Hq f)(ax)$. Hence
\begin{align*}
\ip{h_a(\Hq\LL)}{f}&=-\ip{\LL}{\Hq(h_af)}=-a\ip{\LL}{h_a(\Hq f)}\\
&=-a\ip{h_a\LL}{\Hq f}=a\ip{\Hq(h_a\LL)}{f}.
\end{align*}

(iv) By the first identity in (ii) applied with $f=h_qR$ and the functional
$R\phi u$,
\[
\Hq\big((h_qR)\phi v\big)=\Hq\big((h_qR)\,R\phi u\big)
=R\,\Hq(R\phi u)+(\Hq R)\,R\phi u .
\]
Applying (ii) again, with $f=R$ and the functional $\phi u$, and using
$\Hq(\phi u)=-\psi u$,
\begin{align*}
R\,\Hq(R\phi u)&=R\Big((\hqi R)\,\Hq(\phi u)+\Hq(\hqi R)\,\phi u\Big)\\
&=-(\hqi R)\,\psi\,v+\phi\,\Hq(\hqi R)\,v .
\end{align*}
Adding the two relations gives the result.
\end{proof}

\begin{remark}
By Lemma \ref{lem:rules}(iv), a polynomial (Christoffel) modification $Ru$
of an $\Hq$-semiclassical functional $u$ is $\Hq$-semiclassical as soon as
it is quasi-definite; the condition on the degree of the second
coefficient is checked as in the proof of Proposition \ref{prop:main}.
\end{remark}

%======================================================================
\section{The functional $\Le$ and positive-definiteness}\label{sec2}
%======================================================================

We now fix a normalized functional $\LL$ and a number
$\epsilon\in\mathbb C$, and we turn to the equation
\begin{equation}\label{eq:main}
\Le-\epsilon\,\Hq\Le=\LL .
\end{equation}

\begin{proposition}\label{prop:moments}
Equation \eqref{eq:main} has a unique solution $\Le\in\PP'$. It is
normalized and
\begin{equation}\label{eq:mom}
(\Le)_n=\sum_{\nu=0}^n(-\epsilon)^{n-\nu}\,
\frac{\qn n!}{\qn\nu!}\,(\LL)_\nu ,\qquad n\ge0 .
\end{equation}
Moreover, put $T_\epsilon p:=\sum_{k\ge0}(-\epsilon)^k\Hq^kp$ for
$p\in\PP$ (a finite sum). Then $T_\epsilon p$ is the unique polynomial
solution $y$ of $y+\epsilon\Hq y=p$, it has the same degree and the same
leading coefficient as $p$, and
\begin{equation}\label{eq:T}
\ip{\Le}{p}=\ip{\LL}{T_\epsilon p},\qquad p\in\PP .
\end{equation}
\end{proposition}

\begin{proof}
By \eqref{eq:momH}, equation \eqref{eq:main} is equivalent to
\[
(\Le)_0=(\LL)_0=1,\qquad
(\Le)_n+\epsilon\,\qn n\,(\Le)_{n-1}=(\LL)_n\quad(n\ge1).
\]
This recurrence determines all the moments, so $\Le$ exists and is unique.
Formula \eqref{eq:mom} holds for $n=0$; if it holds for $n-1$, then
\begin{align*}
(\Le)_n&=(\LL)_n-\epsilon\qn n\sum_{\nu=0}^{n-1}(-\epsilon)^{n-1-\nu}
\frac{\qn{n-1}!}{\qn\nu!}(\LL)_\nu\\
&=(\LL)_n+\sum_{\nu=0}^{n-1}(-\epsilon)^{n-\nu}\frac{\qn n!}{\qn\nu!}(\LL)_\nu ,
\end{align*}
which is \eqref{eq:mom}.

The sum defining $T_\epsilon p$ is finite because $\Hq$ lowers the degree.
It is telescoping:
\[
T_\epsilon p+\epsilon\Hq T_\epsilon p
=\sum_{k\ge0}(-\epsilon)^k\Hq^kp-\sum_{k\ge0}(-\epsilon)^{k+1}\Hq^{k+1}p=p .
\]
If $z\neq0$ is a polynomial of degree $d$ with $z+\epsilon\Hq z=0$, the
coefficient of $x^d$ in $z+\epsilon\Hq z$ is the leading coefficient of
$z$, which is absurd; this gives uniqueness. The terms with $k\ge1$ in
$T_\epsilon p$ have degree $<\deg p$, so $T_\epsilon p$ and $p$ have the
same degree and leading coefficient. Finally, for every $\M\in\PP'$ and
$y\in\PP$, \eqref{eq:dual} gives
$\ip{\M-\epsilon\Hq\M}{y}=\ip{\M}{y+\epsilon\Hq y}$. With $\M=\Le$ and
$y=T_\epsilon p$ we obtain
$\ip{\Le}{p}=\ip{\Le-\epsilon\Hq\Le}{T_\epsilon p}=\ip{\LL}{T_\epsilon p}$.
\end{proof}

With $p=x^n$, \eqref{eq:T} gives back \eqref{eq:mom}. By Lemma
\ref{lem:rules}(iii), if $\Le$ solves \eqref{eq:main}, then $h_a\Le$
solves the same equation with data $h_a\LL$ and parameter $a\epsilon$. In
particular the reflection $h_{-1}$ exchanges $\epsilon$ and $-\epsilon$.

\medskip
For the rest of this section $q$ is real and $0<q<1$. One functional plays
a special part in what follows; we introduce it now and study it in
Section \ref{sec4}:
\begin{equation}\label{eq:K}
\ip{\Kq}{f}=\sum_{k\ge0}w_k\,f(x_k),\qquad
x_k=\frac{q^k}{1-q},\qquad w_k=q^k\,(q^{k+1};q)_\infty>0 .
\end{equation}
For the moment we only need the fact that $(\Kq)_n=\qn n!$
(Lemma \ref{lem:K}).

\begin{lemma}\label{lem:delta}
Let $0<q<1$, $y\in\mathbb R$ and $\epsilon\in\mathbb R\setminus\{0\}$, and
let $(\delta_y)_\epsilon$ be the solution of \eqref{eq:main} with
$\LL=\delta_y$. By \eqref{eq:T},
$\ip{(\delta_y)_\epsilon}{p}=(T_\epsilon p)(y)$.
\begin{enumerate}[label=(\roman*)]
\item $(\delta_0)_\epsilon=h_{-\epsilon}\Kq$.
\item If $\epsilon y>0$, put $s=(1-q)y/\epsilon>0$. Then
\begin{equation}\label{eq:delta}
(\delta_y)_\epsilon=\sum_{j\ge0}\frac{s\,q^j}{(-s;q)_{j+1}}\,
\delta_{yq^j}+\frac{1}{(-s;q)_\infty}\,h_{-\epsilon}\Kq .
\end{equation}
\end{enumerate}
Consequently, if $\epsilon y\ge0$, then $(\delta_y)_\epsilon$ is a positive
measure of total mass $1$, and $(T_\epsilon p)(y)\ge0$ for every real
polynomial $p\ge0$ on $\mathbb R$.
\end{lemma}

\begin{proof}
(i) By \eqref{eq:mom} with $(\delta_0)_\nu=\delta_{\nu,0}$, the moments of
$(\delta_0)_\epsilon$ are $(-\epsilon)^n\qn n!$. These are the moments of
$h_{-\epsilon}\Kq$, since $\ip{h_{-\epsilon}\Kq}{x^n}=(-\epsilon)^n(\Kq)_n$.

(ii) Denote by $\rho$ the right-hand side of \eqref{eq:delta}. Put
\[
t_j=yq^j,\qquad a_j=\frac{\epsilon}{(1-q)t_j}=\frac{1}{s\,q^j}\qquad(j\ge0),
\]
\[
d_{-1}=1,\qquad d_j=\prod_{i=0}^j\frac{a_i}{1+a_i}=\frac{1}{(-s;q)_{j+1}}
\qquad(j\ge0).
\] The coefficient of $\delta_{t_j}$ in \eqref{eq:delta} is
$c_j=sq^j/(-s;q)_{j+1}=d_j/a_j$, and
\begin{equation}\label{eq:cd}
c_j(1+a_j)=d_j\,\frac{1+a_j}{a_j}=d_{j-1}.
\end{equation}
All $c_j$ are positive and $c_j\le s\,q^j$, so $\sum_jc_j<\infty$; since
$h_{-\epsilon}\Kq$ is a finite positive measure, $\rho$ is a finite
positive measure with bounded support, and in particular $\rho\in\PP'$.

Let $f\in\PP$. Since $t_{j+1}=qt_j$,
\[
\epsilon(\Hq f)(t_j)=\epsilon\,\frac{f(t_{j+1})-f(t_j)}{(q-1)t_j}
=-a_j\big(f(t_{j+1})-f(t_j)\big),
\]
and therefore, by \eqref{eq:cd},
\[
c_j\,(f+\epsilon\Hq f)(t_j)=c_j(1+a_j)f(t_j)-c_ja_jf(t_{j+1})
=d_{j-1}f(t_j)-d_jf(t_{j+1}).
\]
Summing over $0\le j\le J$, the sum telescopes to $f(y)-d_Jf(t_{J+1})$.
When $J\to\infty$ we have $t_{J+1}\to0$ and $d_J\to d_\infty:=1/(-s;q)_\infty$,
hence
\[
\sum_{j\ge0}c_j\,(f+\epsilon\Hq f)(t_j)=f(y)-d_\infty f(0).
\]
On the other hand, by (i) and \eqref{eq:dual},
$\ip{h_{-\epsilon}\Kq}{f+\epsilon\Hq f}
=\ip{(\delta_0)_\epsilon-\epsilon\Hq(\delta_0)_\epsilon}{f}=f(0)$.
Adding up, $\ip{\rho}{f+\epsilon\Hq f}=f(y)$ for every $f\in\PP$, that is
$\rho-\epsilon\Hq\rho=\delta_y$. By uniqueness
(Proposition \ref{prop:moments}), $\rho=(\delta_y)_\epsilon$.

The total mass of $\rho$ is
$\sum_{j\ge0}(d_{j-1}-d_j)+d_\infty=d_{-1}=1$, since $c_j=d_{j-1}-d_j$ by
\eqref{eq:cd}. The last assertion follows from (i), (ii) and the positivity
of $\Kq$, because $\ip{(\delta_y)_\epsilon}{p}=(T_\epsilon p)(y)$.
\end{proof}

\begin{remark}
For $q=1$ and $\epsilon>0$, the polynomial solution of $y+\epsilon y'=p$
is $y(x)=\epsilon^{-1}\int_{-\infty}^x e^{(t-x)/\epsilon}p(t)\,dt$, which
is non-negative when $p\ge0$; this is the mechanism behind
\cite[Lemma 2]{Sf}. Formula \eqref{eq:delta} is its $q$-analogue. The
lattice part $\{yq^j\}$ replaces the exponential density, and the term
$h_{-\epsilon}\Kq$ carries the mass that the geometric lattice cannot
carry beyond the origin.
\end{remark}

\begin{theorem}\label{thm:main}
Let $0<q<1$ and $\epsilon>0$. If $\LL$ is positive-definite on
$[0,+\infty)$, then $\Le$ is positive-definite.
\end{theorem}

\begin{proof}
Let $p$ be a non-zero real polynomial with $p\ge0$ on $\mathbb R$. For every
$y\ge0$ we have $\epsilon y\ge0$, so Lemma \ref{lem:delta} gives
$(T_\epsilon p)(y)\ge0$. Thus $T_\epsilon p$ is a real polynomial, which is
non-negative on $[0,+\infty)$ and non-zero, since it has the same degree as
$p$. As $\LL$ is positive-definite on $[0,+\infty)$, \eqref{eq:T} yields
$\ip{\Le}{p}=\ip{\LL}{T_\epsilon p}>0$.
\end{proof}

\begin{corollary}\label{cor:refl}
Let $0<q<1$ and $\epsilon<0$. If $\LL$ is positive-definite on
$(-\infty,0]$, then $\Le$ is positive-definite.
\end{corollary}

\begin{proof}
The functional $h_{-1}\LL$ is positive-definite on $[0,+\infty)$ and, by
Lemma \ref{lem:rules}(iii), $h_{-1}\Le$ is the solution of \eqref{eq:main}
with data $h_{-1}\LL$ and parameter $-\epsilon>0$. By Theorem
\ref{thm:main}, $h_{-1}\Le$ is positive-definite, hence so is $\Le$.
\end{proof}

\begin{remark}
In \cite[Theorem 3]{Sf} ($q=1$), every real $\epsilon$ is allowed and $\LL$
may be any positive-definite functional. For $0<q<1$ the sign of $\epsilon$
has to match the half-line that carries $\LL$, because $\Hq$, unlike the
derivative, does not commute with translations. Theorem \ref{thm:K} shows
that, for $\LL=\Kq$ and $\epsilon<0$, positivity actually fails except when
$\epsilon$ lies on the lattice $-q^{\mathbb Z}$.
\end{remark}

%======================================================================
\section{The $\Hq$-semiclassical character}\label{sec3}
%======================================================================

We now ask when $\Le$ is $\Hq$-semiclassical. The first observation is
that $\Le$ always inherits a second-order equation from $\LL$.

\begin{proposition}\label{prop:order2}
Let $\LL$ satisfy $\Hq(\phi\LL)+\psi\LL=0$. Then $\Le$ satisfies the
second-order equation
\begin{equation}\label{eq:order2}
\Hq\Big(-\epsilon\,\Hq\big((h_q\phi)\Le\big)
+\big(\phi+\epsilon\Hq\phi-\epsilon\,h_q\psi\big)\Le\Big)
+\big(\psi+\epsilon\Hq\psi\big)\Le=0 .
\end{equation}
\end{proposition}

\begin{proof}
Let $g\in\PP$. Multiplying \eqref{eq:main} by $g$ and using \eqref{eq:fH},
\[
g\LL=g\Le-\epsilon\,g\,\Hq\Le
=g\Le-\epsilon\,\Hq\big((h_qg)\Le\big)+\epsilon\,(\Hq g)\Le .
\]
With $g=\phi$ this gives
$\phi\LL=-\epsilon\Hq((h_q\phi)\Le)+(\phi+\epsilon\Hq\phi)\Le$, and with
$g=\psi$ it gives
$\psi\LL=-\epsilon\Hq((h_q\psi)\Le)+(\psi+\epsilon\Hq\psi)\Le$.
Substituting these two expressions into $\Hq(\phi\LL)+\psi\LL=0$ and
gathering the terms under $\Hq$ yields \eqref{eq:order2}.
\end{proof}

When $q\to1$, \eqref{eq:order2} reduces to \cite[(3.2)]{Sf}.

\begin{proposition}\label{prop:main}
Let $\epsilon\neq0$ and let $\LL$ and $\M=\Le$ be quasi-definite
functionals satisfying \eqref{eq:main}. The following statements are
equivalent:
\begin{enumerate}[label=(\roman*)]
\item $\M$ is $\Hq$-semiclassical and satisfies
$\Hq(\phi_1\M)+\psi_1\M=0$;
\item there exist non-zero polynomials $E$ and $F$ such that $E\LL=F\M$.
\end{enumerate}
More precisely, (i) implies (ii) with
\begin{equation}\label{eq:EF}
E=\hqi\phi_1,\qquad F=E+\epsilon\big(\psi_1+\Hq E\big),
\end{equation}
and (ii) implies
\begin{equation}\label{eq:Msc}
\Hq\big((h_qE)\,\M\big)+\epsilon^{-1}\big(F-E-\epsilon\Hq E\big)\M=0 .
\end{equation}
Furthermore, if (ii) holds, then $\LL$ is $\Hq$-semiclassical and
\begin{equation}\label{eq:Lsc}
\Hq\big((h_qF)\,E\,\LL\big)+\Psi\,\LL=0,
\end{equation}
where
\[
\Psi=\epsilon^{-1}(\hqi F)(F-E)-E\big(\Hq F+\Hq(\hqi F)\big).
\]
\end{proposition}

\begin{proof}
\emph{(i)$\Rightarrow$(ii).} Put $E=\hqi\phi_1$, so that $\phi_1=h_qE$. By
\eqref{eq:fH} with $f=E$,
\[
E\,\Hq\M=\Hq(\phi_1\M)-(\Hq E)\M=-\big(\psi_1+\Hq E\big)\M .
\]
Multiplying \eqref{eq:main} by $E$ we get
\[
E\LL=E\M-\epsilon E\,\Hq\M=\big(E+\epsilon(\psi_1+\Hq E)\big)\M=F\M .
\]
$E\neq0$ because $\phi_1\neq0$. If $F$ were zero, then $E\LL=0$, which
contradicts the quasi-definiteness of $\LL$. Hence $F\neq0$.

\emph{(ii)$\Rightarrow$(i).} Assume $E\LL=F\M$. Multiplying
\eqref{eq:main} by $E$ and using \eqref{eq:fH} with $f=E$,
\[
F\M=E\LL=E\M-\epsilon E\,\Hq\M
=E\M-\epsilon\Hq\big((h_qE)\M\big)+\epsilon(\Hq E)\M .
\]
Dividing by $-\epsilon$ we obtain \eqref{eq:Msc}. Put
$\Phi=h_qE\neq0$ and $\Psi_1=\epsilon^{-1}(F-E-\epsilon\Hq E)$. Suppose
that $\Psi_1$ is a constant $c$. Applying \eqref{eq:Msc} to the polynomial
$1$ and using $\ip{\Hq(\Phi\M)}{1}=-\ip{\Phi\M}{\Hq1}=0$, we get
$c\,(\M)_0=0$. Since $\M$ is quasi-definite, $(\M)_0\neq0$, so $c=0$. Then
$\Hq(\Phi\M)=0$, hence $\Phi\M=0$ and $\Phi=0$, a contradiction.
Therefore $\deg\Psi_1\ge1$ and, after division by the leading
coefficient of $\Phi$, \eqref{eq:Msc} shows that $\M$ is
$\Hq$-semiclassical.

\emph{The equation of $\LL$.} From \eqref{eq:main},
$\Hq\M=\epsilon^{-1}(\M-\LL)$. By Lemma \ref{lem:rules}(ii) with $f=F$,
\begin{align*}
\Hq(E\LL)=\Hq(F\M)&=(\hqi F)\,\Hq\M+\Hq(\hqi F)\,\M\\
&=\epsilon^{-1}(\hqi F)(\M-\LL)+\Hq(\hqi F)\,\M .
\end{align*}
We multiply by $F$ and use $F\M=E\LL$ in each term containing $\M$:
\[
F\,\Hq(E\LL)=\epsilon^{-1}(\hqi F)(E-F)\LL+\Hq(\hqi F)\,E\LL .
\]
Finally, \eqref{eq:fH} applied with $f=F$ to the functional $E\LL$ gives
$F\,\Hq(E\LL)=\Hq\big((h_qF)E\LL\big)-(\Hq F)E\LL$. Comparing the two
expressions of $F\,\Hq(E\LL)$ we obtain \eqref{eq:Lsc}. The polynomial
$(h_qF)E$ is non-zero, and $\deg\Psi\ge1$ follows from the quasi-definiteness
of $\LL$ exactly as above. Hence $\LL$ is $\Hq$-semiclassical.
\end{proof}

When $q\to1$, \eqref{eq:EF}, \eqref{eq:Msc} and \eqref{eq:Lsc} reduce to
\cite[Proposition 4, (3.5), (3.6)]{Sf}. One difference should be kept in
mind: no translation commutes with $\Hq$, so only dilations are available
to normalize the examples below.

\subsection{Examples}

Let $\mathcal H_q$ be the symmetric normalized functional satisfying
$\Hq\mathcal H_q+2x\mathcal H_q=0$, that is
$(\mathcal H_q)_{2n+1}=0$ and
$(\mathcal H_q)_{2n}=\prod_{k=1}^n\frac{\qn{2k-1}}2$. For
$\alpha\in\mathbb C$, let $\mathcal L_q(\alpha)$ be the normalized
functional satisfying
\begin{equation}\label{eq:lagq}
\Hq\big(qx\,\mathcal L_q(\alpha)\big)+(x-\alpha-1)\,\mathcal L_q(\alpha)=0,
\qquad
(\mathcal L_q(\alpha))_n=\prod_{k=1}^n\big(\alpha+\qn k\big).
\end{equation}
(The moments follow from \eqref{eq:lagq} applied to $x^n$, since
$q\qn n+1=\qn{n+1}$.) We assume that these functionals are quasi-definite.
They are $\Hq$-classical, and they tend to the Hermite functional
$\mathcal H$ and to the Laguerre functional $\mathcal L(\alpha)$ when
$q\to1$. Note that $\mathcal L_q(0)=\Kq$.

\begin{example}[$\Le=\mathcal H_q$]\label{ex:her}
Here $\phi_1=1$ and $\psi_1=2x$, so \eqref{eq:EF} gives $E=1$ and
$F=1+2\epsilon x$. This can also be seen directly: $\Hq\mathcal
H_q=-2x\mathcal H_q$, hence
\[
\LL=\mathcal H_q-\epsilon\Hq\mathcal H_q=(1+2\epsilon x)\,\mathcal H_q .
\]
This is the $q$-analogue of the Hermite example of \cite[p.\,682]{Sf},
where the functional is printed as
$(2\epsilon)^{-1}\big(x+(2\epsilon)^{-1}\big)\mathcal H$; the constant in
front should be $2\epsilon$, which gives $(1+2\epsilon x)\mathcal H$. As a Christoffel transform,
$\LL$ is quasi-definite if and only if the monic orthogonal polynomials
of $\mathcal H_q$ do not vanish at $-(2\epsilon)^{-1}$. Since
$h_qF=1+2\epsilon qx$, $\hqi F=1+2\epsilon q^{-1}x$,
$\Hq F=2\epsilon$ and $\Hq(\hqi F)=2\epsilon q^{-1}$, equation
\eqref{eq:Lsc} reads
\[
\Hq\big((1+2\epsilon qx)\LL\big)
+\Big(2x\big(1+2\epsilon q^{-1}x\big)-2\epsilon\big(1+q^{-1}\big)\Big)\LL=0 .
\]
\end{example}

\begin{example}[$\Le=\mathcal L_q(\alpha)$]\label{ex:lag}
Here $\phi_1=qx$ and $\psi_1=x-\alpha-1$. Then $E=\hqi\phi_1=x$,
$\Hq E=1$ and \eqref{eq:EF} gives $F=(1+\epsilon)x-\epsilon\alpha$, so
that
\begin{equation}\label{eq:37}
x\LL=(1+\epsilon)(x-\tilde a)\,\mathcal L_q(\alpha),\qquad
\tilde a=\frac{\epsilon\alpha}{1+\epsilon}\quad(\epsilon\neq-1).
\end{equation}
This is exactly the $q$-analogue of \cite[(3.7)]{Sf}. From
\eqref{eq:momH}, the moments of
$\LL=\mathcal L_q(\alpha)-\epsilon\Hq\mathcal L_q(\alpha)$ are
\[
(\LL)_n=\big(\alpha+(1+\epsilon)\qn n\big)\prod_{k=1}^{n-1}(\alpha+\qn k),
\qquad n\ge1 .
\]

\emph{Case $\epsilon=-1$.} Then $(\LL)_n=\prod_{k=0}^{n-1}(\alpha+\qn k)$.
Since $\qn k=1+q\qn{k-1}$, we have
$\alpha+\qn{k-1}=q^{-1}\big(q\alpha-1+\qn k\big)$, hence
\[
(\LL)_n=q^{-n}\prod_{k=1}^n\big(q\alpha-1+\qn k\big),
\qquad\text{i.e.}\qquad \LL=\hqi\,\mathcal L_q(q\alpha-1).
\]
This is the $q$-analogue of $\LL=\mathcal L(\alpha-1)$ in \cite{Sf}.
Moreover \eqref{eq:EF} gives $x\LL=\alpha\mathcal L_q(\alpha)$, so that
\begin{equation}\label{eq:lagrel}
\mathcal L_q(\alpha)=\alpha^{-1}x\,\hqi\mathcal L_q(q\alpha-1)\qquad
(\alpha\neq0).
\end{equation}

\emph{Case $\epsilon\neq-1$ and $\alpha=0$: a $q$-Uvarov transform.} Then
$\tilde a=0$ and $x\LL=(1+\epsilon)x\Kq$, so
$x\big(\LL-(1+\epsilon)\Kq\big)=0$ and
\[
\LL=(1+\epsilon)\Kq-\epsilon\,\delta_0 ,
\]
the constant being fixed by $(\LL)_0=1$. Write
$\LL=(1+\epsilon)(\Kq+\mu\delta_0)$ with $\mu=-\epsilon/(1+\epsilon)$. By
\cite{MM}, $\Kq+\mu\delta_0$ is quasi-definite if and only if
$1+\mu K_n(0,0)\neq0$ for every $n\ge0$, where
$K_n(0,0)=\sum_{\nu=0}^n\widehat P_\nu(0)^2/\ip{\Kq}{\widehat P_\nu^2}$ and
$\{\widehat P_n\}$ is the monic orthogonal sequence of $\Kq$. By Lemma
\ref{lem:K}(iii), each term of this sum equals $q^{-\nu}$, so
$K_n(0,0)=\sum_{\nu=0}^nq^{-\nu}=q^{-n}\qn{n+1}$. Since
$\qn{n+1}-q^n=\qn n$, the condition $1+\mu K_n(0,0)\neq0$ is equivalent to
$q^n\neq\epsilon\qn n$. Therefore
\[
\LL\ \text{is quasi-definite}\iff
\epsilon\neq-1\ \text{and}\ \epsilon\neq\frac{q^n}{\qn n}\quad(n\ge1).
\]
When $q\to1$ this is the condition $\epsilon\neq1/n$ of \cite{Sf}. In this
case $E=x$ and $F=(1+\epsilon)x$, so \eqref{eq:Lsc}, divided by
$(1+\epsilon)q$, becomes
\[
\Hq\big(x^2\LL\big)+q^{-2}\,x\,(x-1-q)\,\LL=0 .
\]

\emph{Case $\epsilon\neq-1$ and $\alpha\neq0$.} By \eqref{eq:37} and
\eqref{eq:lagrel},
$x\LL=x\,\alpha^{-1}(1+\epsilon)(x-\tilde a)\,\hqi\mathcal L_q(q\alpha-1)$.
Both functionals in this identity have the same moment of order $0$, so
\[
\LL=\alpha^{-1}(1+\epsilon)(x-\tilde a)\,\hqi\mathcal L_q(q\alpha-1).
\]
It is a Christoffel transform, quasi-definite if and only if the monic
orthogonal polynomials of $\hqi\mathcal L_q(q\alpha-1)$ do not vanish at
$\tilde a$.
\end{example}

%======================================================================
\section{The equation $\Le-\epsilon\Hq\Le=\Kq$}
\label{sec4}
%======================================================================

In this last section the data is the functional $\Kq$ introduced in
\eqref{eq:K}. We begin with its basic properties.

\begin{lemma}\label{lem:K}
Let $0<q<1$ and let $\Kq$ be defined by \eqref{eq:K}.
\begin{enumerate}[label=(\roman*)]
\item $(\Kq)_n=\qn n!$ for $n\ge0$, and $\Kq$ is positive-definite.
\item $\Kq$ is $\Hq$-classical:
$\Hq(qx\,\Kq)+(x-1)\Kq=0$. Moreover $\Hq\Kq=-\Kq+\delta_0$.
\item The monic orthogonal sequence $\{\widehat P_n\}$ of $\Kq$ satisfies
$\widehat P_{n+1}=(x-\beta_n)\widehat P_n-\gamma_n\widehat P_{n-1}$ with
\[
\beta_n=q^n\big(\qn n+\qn{n+1}\big),\qquad
\gamma_{n+1}=q^{2n+1}\qn{n+1}^2 ,
\]
and
\[
\widehat P_n(0)=(-1)^nq^{n(n-1)/2}\qn n!,\qquad
\ip{\Kq}{\widehat P_n^2}=q^{n^2}\big(\qn n!\big)^2 .
\]
\end{enumerate}
\end{lemma}

\begin{proof}
(i) We use the $q$-binomial theorem
\[
\sum_{k\ge0}\frac{z^k}{(q;q)_k}=\frac{1}{(z;q)_\infty}\qquad(|z|<1),
\]
together with $(q^{k+1};q)_\infty=(q;q)_\infty/(q;q)_k$:
\[
\sum_{k\ge0}w_kx_k^n
=\frac{(q;q)_\infty}{(1-q)^n}\sum_{k\ge0}\frac{q^{k(n+1)}}{(q;q)_k}
=\frac{(q;q)_\infty}{(1-q)^n\,(q^{n+1};q)_\infty}
=\frac{(q;q)_n}{(1-q)^n}=\qn n! .
\]
$\Kq$ is a positive measure with infinite support, hence positive-definite.

(ii) Applied to $x^n$, the equation $\Hq(qx\Kq)+(x-1)\Kq=0$ reads
$-q\qn n\qn n!+\qn{n+1}!-\qn n!=0$, which holds because
$\qn{n+1}=1+q\qn n$. By \eqref{eq:momH},
$(\Hq\Kq)_n=-\qn n\qn{n-1}!=-\qn n!$ for $n\ge1$ and $(\Hq\Kq)_0=0$;
these are the moments of $-\Kq+\delta_0$.

(iii) The map $x\mapsto(1-q)x$ sends $\Kq$ to the orthogonality measure
\[
\sum_{k\ge0}\frac{(q;q)_\infty\,q^k}{(q;q)_k}\,\delta_{q^k}
\]
of the little
$q$-Laguerre (Wall) polynomials with parameter $a=1$ \cite[Ch.\,14]{KLS}.
Their recurrence coefficients, rescaled, give $\beta_n$ and $\gamma_n$. The
value of $\widehat P_n(0)$ follows by induction from
$\widehat P_{n+1}(0)=-\beta_n\widehat P_n(0)-\gamma_n\widehat P_{n-1}(0)$:
if $\widehat P_k(0)=(-1)^kq^{k(k-1)/2}\qn k!$ for $k=n-1,n$, then
\begin{align*}
\widehat P_{n+1}(0)&=(-1)^{n+1}q^{n(n-1)/2}\qn n!\;
q^n\big(\qn n+\qn{n+1}-\qn n\big)\\
&=(-1)^{n+1}q^{n(n+1)/2}\qn{n+1}! ,
\end{align*}
where we used
\[
\gamma_n\,q^{(n-1)(n-2)/2}\qn{n-1}!=q^{2n-1+\frac{(n-1)(n-2)}2}\qn n\qn n!
=q^{n+\frac{n(n-1)}2}\qn n\qn n! .
\] Finally
$\ip{\Kq}{\widehat P_n^2}=\gamma_1\cdots\gamma_n
=\prod_{k=1}^nq^{2k-1}\qn k^2=q^{n^2}(\qn n!)^2$.
\end{proof}

Thus $\Kq$ is the $q$-analogue of the Laguerre functional $\mathcal L(0)$,
whose weight is $e^{-x}$. In the notation of \cite{KLS}, \eqref{eq:K} is the
Jackson integral
\[
\int_0^{1/(1-q)}f(x)\,E_q\big(-q(1-q)x\big)\,d_qx ,\qquad
E_q(z)=(-z;q)_\infty .
\]
We now solve
\begin{equation}\label{eq:42}
\Le-\epsilon\Hq\Le=\Kq,\qquad \epsilon\in\mathbb R .
\end{equation}
By \eqref{eq:mom} with $(\Kq)_\nu=\qn\nu!$, exactly as in \cite[(4.3)]{Sf},
\begin{equation}\label{eq:43}
(\Le)_n=\qn n!\sum_{k=0}^n(-\epsilon)^k
=\begin{cases}
\dfrac{1-(-\epsilon)^{n+1}}{1+\epsilon}\,\qn n!,&\epsilon\neq-1,\\[3mm]
(n+1)\,\qn n!,&\epsilon=-1.
\end{cases}
\end{equation}

\begin{theorem}\label{thm:K}
Let $0<q<1$ and $\epsilon\in\mathbb R\setminus\{0,-1\}$. Then
\begin{equation}\label{eq:44}
\Le=\frac{1}{1+\epsilon}\Big(\Kq+\epsilon\,h_{-\epsilon}\Kq\Big),
\end{equation}
that is,
\[
\ip{\Le}{p}=\frac{1}{1+\epsilon}\sum_{k\ge0}w_k\big(p(x_k)
+\epsilon\,p(-\epsilon x_k)\big),\qquad p\in\PP .
\]
Consequently:
\begin{enumerate}[label=(\roman*)]
\item If $\epsilon>0$, $\Le$ is a positive discrete measure supported by
$\{x_k\}\cup\{-\epsilon x_k\}$; it is positive-definite on $\mathbb R$ but
not on $[0,+\infty)$.
\item If $\epsilon=-q^j$ with $j\ge1$, then $\Le=R_j\,\Kq$, where
\[
R_j(x)=\frac{1}{1-q^j}\Big(1-\prod_{i=1}^j
\big(1-q^{i-j}(1-q)x\big)\Big),\qquad\deg R_j=j,
\]
and $R_j(x_m)>0$ for every $m\ge0$.
\item If $\epsilon=-q^{-l}$ with $l\ge1$, then $\Le$ is the positive
discrete measure with masses
\[
\frac{q^m\big((q^{m+1};q)_\infty-(q^{m+l+1};q)_\infty\big)}{1-q^{-l}}
\]
at the points $x_m$, $m\ge-l$, where $(q^{m+1};q)_\infty=0$ for $m<0$. It
is a polynomial modification of degree $l$ of $h_{q^{-l}}\Kq$.
\item If $\epsilon<0$ and $\epsilon\notin-q^{\mathbb Z}$, then $\Le$ is not
positive-definite.
\end{enumerate}
In cases (ii) and (iii), $\Le$ is positive-definite and
$\Hq$-semiclassical. In particular
\[
\mathfrak L_{-q}=x\,\Kq=h_q\mathcal L_q(q^{-1}),
\]
which is the $q$-analogue of $\mathfrak L_{-1}=x\mathcal L(0)=\mathcal L(1)$
in \cite{Sf}.
\end{theorem}

\begin{proof}
\emph{Formula \eqref{eq:44}.} Since
$\ip{h_{-\epsilon}\Kq}{x^n}=(-\epsilon)^n\qn n!$, the moments of the
right-hand side of \eqref{eq:44} are
\[
\frac{1+\epsilon(-\epsilon)^n}{1+\epsilon}\,\qn n!
=\frac{1-(-\epsilon)^{n+1}}{1+\epsilon}\,\qn n!,
\]
which are the moments \eqref{eq:43}.

\emph{(i)} For $\epsilon>0$ all the masses $w_k/(1+\epsilon)$ and
$\epsilon w_k/(1+\epsilon)$ are positive, and the support is infinite.
Hence $\Le$ is positive-definite. Its moments satisfy
$|(\Le)_n|\le C\,R^n$ with $R=\max(1,\epsilon)/(1-q)$, so the moment problem
is determinate and \eqref{eq:44} is the only representing measure. If
$\Le$ were positive-definite on $[0,+\infty)$, it would have a representing
measure supported in $[0,+\infty)$ (Stieltjes moment problem, see
\cite{Ch}); by determinacy this measure would be \eqref{eq:44}, which
charges $(-\infty,0)$. Hence $\Le$ is not positive-definite on
$[0,+\infty)$.

\emph{(ii)} Let $\epsilon=-q^j$, $j\ge1$. Then $-\epsilon x_k=x_{k+j}$ and
$1+\epsilon=1-q^j>0$. Collecting the masses at the point $x_m$, the mass of
$\Le$ at $x_m$ equals
\[
\frac{w_m-q^jw_{m-j}}{1-q^j},\qquad\text{with } w_{m-j}:=0 \text{ for } m<j .
\]
Since
\[
(q^{m-j+1};q)_\infty=(q^{m-j+1};q)_j\,(q^{m+1};q)_\infty,
\]
we have
\[
q^jw_{m-j}=q^m(q^{m-j+1};q)_\infty=w_m\,(q^{m-j+1};q)_j ;
\]
this also holds
for $m<j$, because then $(q^{m-j+1};q)_j$ contains the factor $1-q^0=0$.
Hence the mass at $x_m$ is
\[
w_m\,\frac{1-(q^{m-j+1};q)_j}{1-q^j}=w_m\,R_j(x_m),
\]
because $q^{m-j+i}=q^{i-j}(1-q)x_m$. Therefore $\Le=R_j\Kq$. If $m\ge j$,
every factor $1-q^{m-j+i}$, $1\le i\le j$, lies in $(0,1)$, so
$(q^{m-j+1};q)_j\in(0,1)$; if $m<j$, the factor with $i=j-m$ vanishes and
the product is $0$. In both cases $R_j(x_m)>0$.

\emph{(iii)} Let $\epsilon=-q^{-l}$, $l\ge1$. Then $-\epsilon x_k=x_{k-l}$
and $1+\epsilon=1-q^{-l}<0$. The mass of $\Le$ at $x_m$, $m\ge-l$, is
$(w_m-q^{-l}w_{m+l})/(1-q^{-l})$, with $w_m:=0$ for $m<0$, and
$q^{-l}w_{m+l}=q^m(q^{m+l+1};q)_\infty$, which gives the stated formula. For
$m\ge0$ the numerator is negative because
\[
(q^{m+1};q)_\infty=(q^{m+1};q)_l\,(q^{m+l+1};q)_\infty<(q^{m+l+1};q)_\infty .
\] For $-l\le m<0$ the numerator is
$-q^m(q^{m+l+1};q)_\infty<0$. In both cases the mass is positive. Finally,
the masses of $h_{q^{-l}}\Kq$ at $x_m$ are $q^{m+l}(q^{m+l+1};q)_\infty$,
so the mass of $\Le$ at $x_m$ equals that of $h_{q^{-l}}\Kq$ multiplied by
$q^{-l}\big((q^{m+1};q)_l-1\big)/(1-q^{-l})$. This is the value at $x_m$ of
a polynomial of degree $l$, since $q^{m+i}=q^{i}(1-q)\,x_m$.

\emph{(iv)} Let $\epsilon<0$ with $\epsilon\notin-q^{\mathbb Z}$. Then the
lattices $\{x_k\}$ and $\{-\epsilon x_k\}$ are disjoint, so the masses in
\eqref{eq:44} do not combine. If $-1<\epsilon<0$, the masses
$\epsilon w_k/(1+\epsilon)$ are negative; if $\epsilon<-1$, the masses
$w_k/(1+\epsilon)$ are negative. Assume that $\Le$ is positive-definite.
Then $\Le$ is represented by a positive measure $\sigma$ \cite{Ch}. Since
$|(\Le)_n|\le CR^n$ with $R=\max(1,|\epsilon|)/(1-q)$, we have
$\int x^{2n}\,d\sigma\le CR^{2n}$ for all $n$, and therefore $\sigma$ is
supported in $[-R,R]$. The signed measure \eqref{eq:44} has finite total variation and
is supported in $[-R,R]$, and it coincides with $\sigma$ on polynomials. By
the Weierstrass approximation theorem the two measures are equal, which
contradicts the existence of negative masses.

\emph{Positive-definiteness and semiclassical character in (ii)--(iii).} In
both cases $\Le$ is a positive measure with infinitely many mass points,
hence it is positive-definite. It is a polynomial modification of $\Kq$, or
of $h_{q^{-l}}\Kq$, which are $\Hq$-classical by Lemma \ref{lem:K}(ii) and
Lemma \ref{lem:rules}(iii). By Lemma \ref{lem:rules}(iv) and the remark
following it, $\Le$ is $\Hq$-semiclassical.

\emph{The case $\epsilon=-q$.} By \eqref{eq:43},
$(\mathfrak L_{-q})_n=\qn n!\,\frac{1-q^{n+1}}{1-q}=\qn{n+1}!
=(x\Kq)_n$. Moreover $\qn{k+1}=1+q\qn k=q\big(q^{-1}+\qn k\big)$, so
$\qn{n+1}!=q^n\prod_{k=1}^n(q^{-1}+\qn k)
=\ip{h_q\mathcal L_q(q^{-1})}{x^n}$.
\end{proof}

\begin{proposition}\label{prop:notsc}
Let $0<q<1$ and
$\epsilon\in\mathbb R\setminus\big(\{0\}\cup-q^{\mathbb Z}\big)$. Then
$\Le$ is not $\Hq$-semiclassical.
\end{proposition}

\begin{proof}
By \eqref{eq:44}, $\Le=\sum_t\nu_t\,\delta_t$, where $t$ runs over the
union of the lattices $L_1=\{x_k\}_{k\ge0}$ and
$L_2=\{-\epsilon x_k\}_{k\ge0}$. These lattices are disjoint because
$\epsilon\notin-q^{\mathbb Z}$. The masses are
$\nu_{x_k}=w_k/(1+\epsilon)$ and $\nu_{-\epsilon x_k}=\epsilon
w_k/(1+\epsilon)$; they are all non-zero and $\sum_t|\nu_t|<\infty$.

Assume that $\Hq(\phi\Le)+\psi\Le=0$ for some polynomials $\phi,\psi$.
Put $B(t)=\phi(t)/\big((q-1)t\big)$ and $A(t)=\psi(t)+B(t)$. For every
$f\in\PP$,
\begin{equation}\label{eq:pf1}
0=\ip{\Le}{\psi f-\phi\,\Hq f}
=\sum_t\nu_t\Big(\psi(t)f(t)+B(t)\big(f(t)-f(qt)\big)\Big).
\end{equation}
Here $B$ is well defined on the support, which does not contain $0$, and
$qt$ belongs to the support whenever $t$ does.
Let $I$ be a compact interval containing $0$ and the support of $\Le$.
Since $|f(t)-f(qt)|\le(1-q)|t|\,\|f'\|_{\infty,I}$, the right-hand side of
\eqref{eq:pf1} is a continuous linear form on $C^1(I)$ for the norm
$\|f\|_{\infty,I}+\|f'\|_{\infty,I}$. Polynomials are dense in $C^1(I)$ for
this norm, so \eqref{eq:pf1} holds for every $f\in C^1(I)$.

Every point $s\neq0$ of the support is isolated in the support. Choose
$f\in C^1(I)$ with $f(s)=1$ and $f=0$ at all the other points of the
support. In \eqref{eq:pf1}, the only non-zero contributions then come from
$t=s$ and, if $s/q$ belongs to the support, from $t=s/q$. We get
\begin{equation}\label{eq:pf2}
\nu_s\,A(s)=\nu_{s/q}\,B(s/q).
\end{equation}

Take $s=x_{k+1}\in L_1$, so that $s/q=x_k$. Since
$w_{k+1}/w_k=q/(1-q^{k+1})$ and $q^{k+1}=q(1-q)x_k$, \eqref{eq:pf2} gives
\[
q\,A(qx_k)=\big(1-q(1-q)x_k\big)\,B(x_k)\qquad(k\ge0).
\]
Take $s=-\epsilon x_{k+1}\in L_2$, and write $t_k=-\epsilon x_k$. The same
ratio $w_{k+1}/w_k$ now equals $q/\big(1+q(1-q)\epsilon^{-1}t_k\big)$, and
\eqref{eq:pf2} gives
\[
q\,A(qt_k)=\big(1+q(1-q)\epsilon^{-1}t_k\big)\,B(t_k)\qquad(k\ge0).
\]
Both relations are identities between rational functions that hold at
infinitely many points, so they hold identically. Subtracting them, we get
$-q(1-q)\big(1+\epsilon^{-1}\big)\,t\,B(t)\equiv0$. Since $\epsilon\neq-1$,
$B\equiv0$, i.e.\ $\phi=0$. Then $A=\psi$ and the first identity gives
$\psi\equiv0$. Hence no non-trivial relation exists, and $\Le$ is not
$\Hq$-semiclassical.
\end{proof}

The proof gives slightly more than the statement: for these values of
$\epsilon$ there is no pair $(\phi,\psi)\neq(0,0)$ of polynomials with
$\Hq(\phi\Le)+\psi\Le=0$, whether $\Le$ is quasi-definite or not.

It remains to look at $\epsilon=-1$, which is excluded from Theorem
\ref{thm:K}. By \eqref{eq:43}, $(\mathfrak L_{-1})_n=(n+1)\qn n!$, so that
\begin{equation}\label{eq:Lm1}
\ip{\mathfrak L_{-1}}{p}=\sum_{k\ge0}w_k\,(xp)'(x_k),\qquad p\in\PP .
\end{equation}
This functional involves derivatives at the mass points, and it behaves
like the non-lattice values.

\begin{proposition}\label{prop:minus1}
Let $0<q<1$. The functional $\mathfrak L_{-1}$ is not positive-definite,
and there is no pair $(\phi,\psi)\neq(0,0)$ of polynomials such that
$\Hq(\phi\mathfrak L_{-1})+\psi\mathfrak L_{-1}=0$. In particular
$\mathfrak L_{-1}$ is not $\Hq$-semiclassical.
\end{proposition}

\begin{proof}
Let $R=1/(1-q)=x_0$ and $I=[-R,R]$. By \eqref{eq:Lm1}, the form
$\Lambda(f)=\sum_kw_k\big(f(x_k)+x_kf'(x_k)\big)$ is continuous on $C^1(I)$
for the norm $\|f\|_{\infty,I}+\|f'\|_{\infty,I}$, and it coincides with
$\mathfrak L_{-1}$ on polynomials.

\emph{Positivity.} Assume that $\mathfrak L_{-1}$ is positive-definite, and
let $\sigma$ be a representing positive measure. Since
$(\mathfrak L_{-1})_{2n}\le(2n+1)R^{2n}$, $\sigma$ is supported in $I$. The
forms $\Lambda$ and $f\mapsto\int f\,d\sigma$ are continuous on $C^1(I)$ and
agree on polynomials, which are dense in $C^1(I)$; hence they agree on
$C^1(I)$. Let $\chi$ be a $C^1$ function on $\mathbb R$ with values in
$[0,1]$, equal to $1$ near $0$ and vanishing outside $[-1,1]$, and put
$f_\delta(x)=(x-x_0)\,\chi\big((x-x_0)/\delta\big)$, where $\delta>0$ is so
small that $x_1<x_0-\delta$. Then $f_\delta$ and $f_\delta'$ vanish at all
the points $x_k$, $k\ge1$, while $f_\delta(x_0)=0$ and $f_\delta'(x_0)=1$.
Thus $\Lambda(f_\delta)=w_0x_0>0$ for every such $\delta$, whereas
$\big|\int f_\delta\,d\sigma\big|\le\delta\,\sigma(I)$ tends to $0$ with
$\delta$. This is a contradiction.

\emph{Pearson equation.} Assume that
$\Hq(\phi\mathfrak L_{-1})+\psi\mathfrak L_{-1}=0$. Put
$H=\phi/(q-1)$ and $G(t)=t\psi(t)+H(t)$. For $f\in\PP$ the polynomial
$h=\psi f-\phi\,\Hq f$ satisfies $t\,h(t)=t\psi(t)f(t)+H(t)\big(f(t)-f(qt)\big)$,
and \eqref{eq:Lm1} gives
\begin{equation}\label{eq:pf3}
0=\sum_{k\ge0}w_k\Big((t\psi f)'+H'\,\big(f-f(q\,\cdot)\big)
+H\,\big(f'-q\,f'(q\,\cdot)\big)\Big)(x_k).
\end{equation}
The right-hand side is continuous on $C^1(I)$, so \eqref{eq:pf3} holds for
every $f\in C^1(I)$. Fix $m\ge1$ and recall that $qx_{m-1}=x_m$. Choosing
$f\in C^1(I)$ with $f'(x_m)=1$, $f(x_m)=0$ and $f=f'=0$ at the other points
$x_k$, and then $f$ with $f(x_m)=1$, $f'(x_m)=0$ and $f=f'=0$ at the other
points $x_k$, we obtain
\[
w_m\,G(x_m)=q\,w_{m-1}H(x_{m-1}),\qquad
w_m\,G'(x_m)=w_{m-1}H'(x_{m-1}).
\]
Since $w_m/w_{m-1}=q/\big(1-q(1-q)x_{m-1}\big)$, this means that, for
$t=x_{m-1}$,
\[
G(qt)=\big(1-q(1-q)t\big)H(t),\qquad
q\,G'(qt)=\big(1-q(1-q)t\big)H'(t).
\]
These are polynomial identities which hold at infinitely many points, hence
identically. Differentiating the first one and comparing with the second
gives $q(1-q)H\equiv0$. Therefore $\phi=0$, then $G=0$, and finally
$\psi=0$.
\end{proof}

\begin{remark}
\begin{enumerate}[label=(\alph*)]
\item When $q\to1$, the lattice $-q^{\mathbb Z}$ fills $(-\infty,0)$, the
discrete measures \eqref{eq:44} become the densities of \cite[Theorem 6]{Sf},
and positive-definiteness for every real $\epsilon$ is recovered. For
$0<q<1$, positive-definiteness and the semiclassical character of $\Le$
for $\epsilon<0$ are arithmetic properties of $\epsilon$.
\item By Proposition \ref{prop:order2} with $\phi=qx$ and $\psi=x-1$, the
functional $\Le$ satisfies a second-order $q$-difference equation for every
$\epsilon$. Hence, for $\epsilon>0$, and more generally for
$\epsilon\notin\{0\}\cup-q^{\mathbb Z}$ whenever $\Le$ is quasi-definite,
$\Le$ satisfies a second-order equation of the form \eqref{eq:order2} but
no first-order equation \eqref{eq:Pearson}; in this sense it is a strict
second-order functional. This is the analogue of
\cite[Proposition 5]{Sf}.
\item For $\epsilon=-1$ the lack of positivity can already be seen on the
first Hankel determinants. One has
$(\mathfrak L_{-1})_2-(\mathfrak L_{-1})_1^2=3q-1$, which is negative for
$q<1/3$ and zero for $q=1/3$; for $q=1/2$ this determinant is positive, but
$\det\big((\mathfrak L_{-1})_{i+j}\big)_{0\le i,j\le2}=-9/128$.
\end{enumerate}
\end{remark}

\section*{Acknowledgements}

The author extends his appreciation to the Deanship of Postgraduate Studies and Scientific Research at Majmaah University for supporting this research work.

\section*{Declarations}

\noindent\textbf{Author contributions.} The author is the sole contributor to this work.

\smallskip
\noindent\textbf{Conflict of interest.} The author declares that he has no competing interests.

\smallskip
\noindent\textbf{Data availability.} No data was used for the research described in the article.

\end{document}